\documentclass{article}

\newtheorem{theorem}{Theorem}[section]
\newtheorem{corollary}[theorem]{Corollary}

\newtheorem{lemma}[theorem]{Lemma}
\newtheorem{proposition}[theorem]{Proposition}

\evensidemargin\oddsidemargin
\begin{document}

\author{Vadim E. Levit and Eugen Mandrescu \\
Department of Computer Science\\
Holon Academic Institute of Technology\\
52 Golomb Str., P.O. Box 305\\
Holon 58102, ISRAEL\\
\{levitv, eugen\_m\}@barley.cteh.ac.il}
\title{On the Structure of $\alpha $-Stable Graphs}
\date{}
\maketitle

\begin{abstract}
The stability number $\alpha (G)$ of a graph $G$ is the cardinality of a
stability system of $G$ (that is of a stable set of maximum size). A graph
is $\alpha $-stable if its stability number remains the same upon both the
deletion and the addition of any edge. Trying to generalize some stable
trees properties, we show that there does not exist any $\alpha $-stable
chordal graph, and we prove that: if $G$ is a connected bipartite graph,
then the following assertions are equivalent: $G$ is $\alpha $-stable; $G$
can be written as a vertex disjoint union of connected bipartite graphs,
each of them having exactly two stability systems covering its vertex set; $G
$ has perfect matchings and $\cap \{M:M$ is a perfect matching of $%
G\}=\emptyset $; from each vertex of $G$ are issuing at least two edges,
contained in some perfect matchings of $G$; any vertex of $G$ belongs to a
cycle, whose edges are alternately in and not in a perfect matching of $G$; $%
\cap \{S:S$ is a stability system of $G\}=\emptyset =\cap \{M:M$ is a
maximum matching of $G\}$.
\end{abstract}
\section{Introduction}

Throughout this paper $G=(V,E)$ is a simple (i.e., a finite, undirected,
loopless and without multiple edges) graph with vertex set $V=V(G)$ and edge
set $E=E(G)$; the number $|V|$ defines its \textit{order}. If $X\subset V$,
then $G[X]$ denotes the subgraph of $G$ spanned by $X$. By $G-W$ we mean the
subgraph $G[V-W]$, if $W\subset V(G)$. By $G-F$ we denote the partial graph
of $G$ obtained by deleting all the edges of $F$, whenever $F\subset E(G)$,
and we use $G-e$, if $W=\{e\}$. If $A,B$ are disjoint subsets of $V(G)$,
then $(A,B)$ stands for the set $\{e=ab:a\in A,b\in B,e\in E\}$. The \textit{%
neighborhood} of a vertex $v\in V$, denoted by $N(v)$, is the set of
vertices adjacent to $v$. For any $X\subset V(G)$, we denote $N(X)=\cup
\{N(x):x\in X\}$. A set of mutually nonadjacent vertices from $V(G)$ is
called \textit{stable} in $G$. A stable set of maximum size will be referred
to as a \textit{stability system} of $G$. The \textit{stability number} of $%
G $, denoted by $\alpha (G)$, is the cardinality of a stability system of $G$%
. A \textit{matching} is a set of non-incident edges of $G$; a matching of
maximum cardinality $\mu (G)$ is a \textit{maximum matching}, and a \textit{%
perfect matching} is a matching covering all the vertices of $G$. A cycle $C$
is called \textit{alternating} if its edges are alternately in and not in a
perfect matching of $G$. A \textit{bipartite graph} is a triple $G=(A,B,E)$,
where $\{A,B\}$ is the bipartition of $V(G)$ to its color classes and $%
E=E(G) $; if $|A|=|B|$, then $G$ is \textit{balanced bipartite}. A set $%
D\subset V$ is called $n$\textit{-dominating} (for $n\geq 1$) in $G$ if $%
|N(v)\cap D|\geq n$, for any vertex $v\in V-D$, \cite{FJ84}. By $C_{n}$, $%
K_{n}$, $P_{n}$ we denote the chordless cycle on $n\geq 4$ vertices, the
complete graph on $n\geq 1$ vertices, and respectively the chordless path on 
$n\geq 3$ vertices.

Based on some results of Haynes et al. from \cite{HLBD90}, Gunther et al.
characterize the trees whose stability numbers are not changing under single
edge addition or deletion (see \cite{GHR93}). In this paper we deal with
chordal graphs (i.e., graphs having no induced $C_{n}$, for $n\geq 4$; see 
\cite{G80}) and bipartite graphs, with the property that their stability
number remains the same upon the deletion and / or the addition of any edge.

In \cite{LM97} the class of $\alpha $-stable graphs, namely, the graphs
whose stability number is insensitive to both deletion and addition of any
edge, has been considered. Here we describe the structure of the $\alpha $%
-stable bipartite graphs.

Some authors have considered similar problems related to adding or deleting
edges and vertices in graphs in connection with various graph parameters
(see \cite{BHNS83}, \cite{BCD88}, \cite{DB88}, \cite{F93}, \cite{HLBD90}, 
\cite{HOS96}, \cite{M96}, \cite{SB83}, \cite{S90}, \cite{TV98}).

\section{$\alpha ^{-}$--stable graphs}

A graph $G$ is called $\alpha ^{-}$\textit{--stable} if $\alpha (G-e)=\alpha
(G)$, for any $e\in E(G)$, (see \cite{GHR93}).

\begin{proposition}
\label{prop2.1}Any graph, which has an $\alpha ^{-}$--stable partial graph
with the same stability number, is an $\alpha ^{-}$--stable graph.
\end{proposition}

\setlength {\parindent}{0.0cm}\textbf{Proof. }Let $H$ be an $\alpha ^{-}$%
--stable partial graph of $G$ with $\alpha (H)=\alpha (G)$. Then, for any $%
e\in E(G)$, we have $\alpha (G)\leq \alpha (G-e)\leq \alpha (H-e)=$ $\alpha
(H)=\alpha (G)$, and this clearly implies the $\alpha ^{-}$--stability of $G$%
. \rule{2mm}{2mm}\setlength
{\parindent}{3.45ex}\newline

Haynes et al. characterized the $\alpha ^{-}$--stable graphs as follows:

\begin{theorem}
\label{th2.2}\cite{HLBD90} A graph $G$ is $\alpha ^{-}$-stable if and only
if each of its stability systems is a $2$-dominating set in $G$.
\end{theorem}

If $G$ is a disconnected graph and $H_{1},...,H_{p}$ are its components,
then $S$ is a stability system of $G$ if and only if $S\cap V(H_{k})$ is a
stability system in $H_{k}$ , for all $k\in \{1,...,p\}$. Consequently, a
disconnected graph is $\alpha ^{-}$-stable if and only if each of its
components is $\alpha ^{-}$-stable.

A graph $G$ having a unique stability system $S$ is called a unique
independence graph; if $V(G)-S$ is also stable, then $G$ is a strong unique
independence graph, (see \cite{HS85}). In this last case, it is clear that $%
G $ is bipartite and has $\{S,V(G)-S\}$ as its bipartition, i.e., its larger
color class equals its unique stability system.

\begin{lemma}
\label{lem2.3}Any unique independence graph is $\alpha ^{-}$-stable.
\end{lemma}

\setlength {\parindent}{0.0cm}\textbf{Proof. }Let $G$ be a graph which has $%
S $ as its unique stability system; it suffices to show that $S$ is also $2$%
-dominating in $G$. Suppose, on the contrary, that $S$ is not such a set.
Then there must exist at least a vertex $v\in V(G)-S$, so that $|N(v)\cap
S|\leq 1$. Now, either: $\mathit{(a)}$ $N(v)\cap S$ $=$ $\emptyset $; then $%
S\cup \{v\}$ is stable in $G$ contradicting the fact that $S$ is a stability
system, or $\mathit{(b)}$ $N(v)\cap S=\{w\}$ and hence $W=S\cup \{v\}-\{w\}$
is a stability system in $G$, again a contradiction, since $S$ is unique.
Consequently, $S$ is a $2$-dominating set in $G$. Hence, by Theorem \ref
{th2.2}, we get that $G$ is $\alpha ^{-}$--stable. \rule{2mm}{2mm}%
\setlength
{\parindent}{3.45ex}\newline

If a graph has a stability system, which is also a $2$-dominating set, then
it is not necessarily unique. However, for graphs in which every even cycle
possesses a chord (so also for chordal graphs), Siemes et al. showed:

\begin{theorem}
\label{th2.4}\cite{STV94} Let $G$ be a graph in which every even cycle has a
chord. If $S$ is a stable set in $G$, then the following statements are
equivalent:

$\mathit{(i)}$ $S$ is the unique stability system of $G$;

$\mathit{(ii)}$ $S$ is a $2$-dominating set of $G$.
\end{theorem}

Combining Lemma 2.3 and Theorem 2.4, we obtain:

\begin{theorem}
\label{th2.5}\cite{LM97} For any chordal graph $G$, the following are
equivalent:

$\mathit{(i)}$ $G$ is $\alpha ^{-}$-stable;

$\mathit{(ii)}$ $G$ has a stability system, which is $2$-dominating;

$\mathit{(iii)}$ $G$ has a unique stability system.
\end{theorem}

In particular, Theorem \ref{th2.5} holds for trees, as Gunther et al. show
in \cite{GHR93}. It is easy to see that this result can not be generalized
to any bipartite graph (e.g., for $n\geq 2$, the complete bipartite graph $%
K_{n,n}$ is $\alpha ^{-}$--stable and has two stability systems). As a
consequence of Proposition \ref{prop2.1}, we get the following result:

\begin{corollary}
\label{cor2.6}Any bipartite graph, which has an $\alpha ^{-}$--stable
spanning tree with the same stability number, is an $\alpha ^{-}$--stable
graph.
\end{corollary}

The bipartite graph consisting of a cycle on $6$ vertices with only one
chord is a counterexample to the converse of the above result. In \cite{HS85}%
, Hopkins and Staton characterize the strong unique independence graphs,
dealing also with the special class of trees. In the sequel, we give some
extensions of their results. Recall that a vertex of degree one is called 
\textit{pendant}.

\begin{lemma}
\label{lem2.7}Any pendant vertex of a graph is contained in some of its
stability systems.
\end{lemma}

\setlength {\parindent}{0.0cm}\textbf{Proof. }Let $v$ be a pendant vertex of
graph $G$ and $S$ be a stability system in $G$. If $v\notin S$ then $S\cap
N(v)=\{w\}$, where $N(v)=\{w\}$ in $G$, and consequently the set $S\cup
\{v\}-\{w\}$ is a stability system in $G$ containing $v$. \rule{2mm}{2mm}%
\setlength {\parindent}{3.45ex}

\begin{corollary}
\label{cor2.8}All pendant vertices of a unique independence graph $G$ are
contained in its unique stability system.
\end{corollary}

\begin{proposition}
\label{prop2.9}If $T$ is a tree of order at least $3$, then the following
conditions are equivalent:

$\mathit{(i)}$ $T$ is a strong unique independence tree;

$\mathit{(ii)}$ $T$ is $\alpha ^{-}$--stable and all its pendant vertices
belong to the larger color class;

$\mathit{(iii)}$ the distance between any two pendant vertices of $T$ is
even.
\end{proposition}

\setlength {\parindent}{0.0cm}\textbf{Proof. }$\mathit{(i)}$ $\Rightarrow $\ 
$\mathit{(ii)}$ Let $S$ be the unique stability system of $T$. Then its
bipartition is $\{S,V-S\}$ and, by above Corollary \ref{cor2.8}, all its
pendant vertices are contained in $S$, i.e., they belong to the larger color
class of $T$.

$\mathit{(ii)}$ $\Rightarrow $\ $\mathit{(iii)}$ It is true, since all
pendant vertices belong to the same color class.

$\mathit{(iii)}$ $\Leftrightarrow $ $\mathit{(i)}$ It has been proved in 
\cite{HS85}. \rule{2mm}{2mm}\setlength {\parindent}{3.45ex}

\begin{lemma}
\label{lem2.10}If $H$ is a partial graph of $G$ such that a stable set $S$
of $G$ is a stability system of $H$, then $S$ is a stability system of $G$
as well.
\end{lemma}

\setlength {\parindent}{0.0cm}\textbf{Proof. }It is clear that $\alpha
(G)\leq \alpha (H)=|S|\leq \alpha (G)$, and this implies the conclusion.\rule%
{2mm}{2mm}\setlength {\parindent}{3.45ex}

\begin{lemma}
\label{lem2.11}If $H$ is a connected partial graph of a bipartite graph $G$,
then $G$ and $H$ have the same bipartition.
\end{lemma}

\setlength {\parindent}{0.0cm}\textbf{Proof. }Let $\{X,Y\}$ be the
bipartition of $H,\{A,B\}$ be the bipartition of $G$ and suppose they are
different. Let denote $A_{1}$ $=$ $A\cup X,A_{2}$ $=$ $A\cup Y,B_{1}$ $=$ $%
B\cup X$ and $B_{2}=B\cup Y$. Since both $G$ and $H$ are connected, at least
one of the sets $(A_{1},A_{2}),(A_{1},B_{2}),(A_{2},B_{2}),(B_{1},B_{2})$
must be non-empty, contradicting the fact that $A,B,X,Y$ are all stable.
Therefore, $G$ and $H$ have the same bipartition. \rule{2mm}{2mm}%
\setlength
{\parindent}{3.45ex}

\begin{proposition}
\label{prop2.12}If $G$ is connected, then the following are equivalent:

$\mathit{(i)}$ $G$ is a strong unique independence graph;

$\mathit{(ii)}$ $G$ is bipartite and has a strong unique independence
spanning tree;

$\mathit{(iii)}$ $G$ is bipartite and has an $\alpha ^{-}$--stable spanning
tree with all its pendant vertices contained in the larger color class;

$\mathit{(iv)}$ $G$ is an $\alpha ^{-}$--stable bipartite graph and one of
its color classes equals its unique stability system.
\end{proposition}

\setlength {\parindent}{0.0cm}\textbf{Proof. }$\mathit{(i)}$ $%
\Leftrightarrow $ $\mathit{(ii)}$ It has been proved in \cite{HS85}.%
\setlength {\parindent}{3.45ex}

$\mathit{(ii)}$ $\Rightarrow $\ $\mathit{(iii)}$ It is true according to
Proposition \ref{prop2.9}.

$\mathit{(iii)}$ $\Rightarrow $\ $\mathit{(iv)}$ Let $T$ be an $\alpha ^{-}$%
--stable spanning tree of $G$, with all its pendant vertices contained in
the larger color class. According to Lemma \ref{lem2.11}, $G$ and $T$ have
the same bipartition, say $\{A,B\}$, and if $A$ contains all pendant
vertices of $T$, then Proposition \ref{prop2.9} ensures that $A$ is its
unique stability system. Hence, by Lemma \ref{lem2.10}, $A$ is a stability
system for $G$ itself. By Corollary \ref{cor2.6}, $G$ is $\alpha ^{-}$%
--stable. If $S\neq A$ is another stability system of $G$, then $S$ is a
stability system of $T$ as well, contradicting the fact that $T$ has a
unique stability system.

$\mathit{(iv)}$ $\Rightarrow $\ $\mathit{(i)}$ It is clear. \rule{2mm}{2mm}%
\setlength {\parindent}{3.45ex}\newline

In order to obtain a characterization for $\alpha ^{-}$-stable bipartite
graphs, we recall the following classical theorem, due to K\"{o}nig.

\begin{theorem}
\label{th2.13}\cite{K31} If $G$ is a bipartite graph of order $n$, then $%
\alpha (G)+\mu (G)=n$.
\end{theorem}

\begin{theorem}
\label{th2.14}If $G$ is a connected bipartite graph, then the following are
equivalent:

$\mathit{(i)}$ $G$ is $\alpha ^{-}$--stable;

$\mathit{(ii)}$ $\mu (G-e)=\mu (G)$, for any $e\in E(G)$;

$\mathit{(iii)}$ $\cap \{M:M$ is a maximum matching of $G\}=\emptyset $.
\end{theorem}

\setlength {\parindent}{0.0cm}\textbf{Proof. }Since $G$ and $G-e$, for any
edge $e\in E(G)$, are bipartite and have the same number of vertices,
Theorem \ref{th2.13} ensures that $\mathit{(i)}$ and $\mathit{(ii)}$ are
equivalent.\setlength {\parindent}{3.45ex}

$\mathit{(i)}$ $\Rightarrow $\ $\mathit{(iii)}$ Let $G$ be $\alpha ^{-}$%
--stable and suppose, on the contrary, that there exists an edge $e$
contained in $\cap \{M:M$ is a maximum matching of $G\}$. Then we have $%
\alpha (G-e)=\alpha (G),\mu (G-e)=\mu (G)-1$, and henceforth $\alpha
(G-e)+\mu (G-e)<|V(G-e)|=|V(G)|$, in contradiction with Theorem \ref{th2.13}.

The converse, $\mathit{(iii)}$ $\Rightarrow $\ $\mathit{(i)}$, can be proven
in a similar way. \rule{2mm}{2mm}

\section{$\alpha ^{+}$-stable graphs}

A graph $G$ is called $\alpha ^{+}$\textit{-stable} if $\alpha (G+e)=\alpha
(G)$, for any $e\in E(\overline{G})$, where $\overline{G}$ is the complement
of $G$ (see \cite{GHR93}). Haynes et al. characterized the $\alpha ^{+}$%
-stable graphs as follows:

\begin{theorem}
\label{th3.1}\cite{HLBD90} A graph $G$ is $\alpha ^{+}$-stable if and only
if no pair of vertices is contained in all its stability systems.
\end{theorem}

An $\alpha ^{+}$-stable graph need not necessarily be connected. If $G$ is a
disconnected graph and $H_{1},...,H_{p}$ are its components, then it is not
difficult to prove that $G$ is $\alpha ^{+}$-stable if and only if each $%
H_{i},i\in \{1,...,p\}$, is $\alpha ^{+}$-stable and at most one of $%
H_{i},i\in \{1,...,p\}$, has $|\cap \{S_{i}:S_{i}$ is a stability system in $%
H_{i}\}|=1$.

Applying Theorem \ref{th3.1} for trees, Gunther et al. proved:

\begin{theorem}
\label{th3.2}\cite{GHR93} For a tree $T$ of order at least two, the
following conditions are equivalent:

$\mathit{(i)}$ $T$ is $\alpha ^{+}$-stable;

$\mathit{(ii)}$ $T$ has two disjoint stability systems that partition its
vertex set;

$\mathit{(iii)}$ $T$ has a perfect matching.
\end{theorem}

\begin{corollary}
\label{cor3.3} For a tree $T$ of order at least two, the next assertions are
equivalent:

$\mathit{(i)}$ $T$ is an $\alpha ^{+}$-stable path;

$\mathit{(ii)}$ $T$ is isomorphic to $P_{2n}$ , for $n\geq 1$;

$\mathit{(iii)}$ the distance between any two pendant vertices of $T$ is odd.
\end{corollary}

\setlength {\parindent}{0.0cm}\textbf{Proof. }The implications $\mathit{(i)}$
$\Leftrightarrow $ $\mathit{(ii)}$ $\Rightarrow $\ $\mathit{(iii)}$ are
obvious.\setlength {\parindent}{3.45ex}

$\mathit{(iii)}$ $\Rightarrow $\ $\mathit{(ii)}$ Since $T$ is a tree, it has
at least two pendant vertices, say $a,b$, which are endpoints of a $P_{2n}$,
for $n\geq $ $1$. If $T\neq P_{2n}$, then it has a third pendant vertex $v$,
and this yields a contradiction: the distance between $v$ and $a$ or between 
$v$ and $b$ is even. \rule{2mm}{2mm}

\begin{lemma}
\label{lem3.4} If $G$ is a connected graph having two stability systems,
which partition its vertex set, then $G$ is bipartite, has a unique
bipartition, up to an isomorphism, and is $\alpha ^{+}$-stable.
\end{lemma}

\setlength {\parindent}{0.0cm}\textbf{Proof. }Clearly, $G$ is bipartite and
the two disjoint stability systems $S_{1}$ and $S_{2}$ generate a
bipartition for its vertex set. Suppose, on the contrary, that $\{A,B\}$ is
another bipartition of $V(G)$, such that $A,B$ are stable in $G$. Then, the
sets $A_{i}=S_{i}\cap A$ and $B_{i}=S_{i}\cap B$, for $i=1,2$, are all
non-empty and clearly we have $S_{i}=A_{i}\cup B_{i},i=1,2$. Therefore,
there is no vertex in $A_{1}\cup B_{2}$, which is adjacent to a vertex from $%
A_{2}\cup B_{1}$, a contradiction, because $G$ is a connected graph. By
Theorem \ref{th3.1}, $G$ is also $\alpha ^{+}$-stable. \rule{2mm}{2mm}%
\setlength {\parindent}{3.45ex}

\begin{lemma}
\label{lem3.5} Any connected bipartite graph has a spanning tree with the
same stability number.
\end{lemma}

\setlength {\parindent}{0.0cm}\textbf{Proof. }Let $G$ be a connected and
bipartite graph. Then $G$ is also perfect and therefore it admits a clique
cover with exactly $\alpha (G)$ cliques (since, by definition of
perfectness, \cite{B73}, any of its induced subgraphs, including itself, has
such a particular cover). Let $H$ be the partial subgraph of $G$ generated
by this clique cover. Evidently, $\alpha (H)=\alpha (G)$ holds. Since $G$ is
bipartite, $H$ consists of a disjoint union of cliques isomorphic to $K_{2}$
and $K_{1}$. Now, add edges to $H$ from $E(G)-E(H)$, such that the new graph 
$T$ is without cycles and connected. This is possible since $G$ was
connected and adding an edge that closes a cycle is clearly redundant in
this rebuilding connectedness process. Clearly, $T$ is a tree and $\alpha
(G)\leq \alpha (T)\leq \alpha (H)$, since $T$ was obtained from $H$ by
adding edges of $G$. Therefore, $T$ is a spanning tree of $G$ with $\alpha
(G)=\alpha (T)$.\rule{2mm}{2mm}\setlength {\parindent}{3.45ex}\newline

This lemma permits us to establish a key relationship between a bipartite
graph and some of its spanning trees.

\begin{proposition}
\label{prop3.6} A connected bipartite graph is $\alpha ^{+}$-stable if and
only if it has an $\alpha ^{+}$-stable spanning tree.
\end{proposition}

\setlength {\parindent}{0.0cm}\textbf{Proof. }Let $G$ be an $\alpha ^{+}$%
-stable connected bipartite graph. By Lemma \ref{lem3.5}, $G$ possesses a
spanning tree $T$ with the same stability number. Now, if $e\in E(G)-E(T)$,
we have $\alpha (T)\geq \alpha (T+e)\geq $ $\alpha (G)$, and for $e\notin
E(G)$, we obtain $\alpha (T)\geq \alpha (T+e)\geq \alpha (G+e)=$ $\alpha (G)$%
. Hence, we get $\alpha (T)=\alpha (T+e)$, for any $e\notin E(T)$, and
consequently, $T$ is $\alpha ^{+}$-stable.\setlength {\parindent}{3.45ex}

Conversely, let $T$ be an $\alpha ^{+}$-stable spanning tree of $G$ and $%
\{A,B\}$ be the bipartition of $V(G)$. Clearly, $\{A,B\}$ is also
bipartition for $V(T)$, because $T$ is a tree obtained by removing edges
from $(A,B)$. Since $T$ is $\alpha ^{+}$-stable, by Theorem \ref{th3.2} and
Lemma \ref{lem3.4}, $T$ has two stability systems, which generate the unique
bipartition of its vertex set. Hence $A,B$ must be the two stability systems
of $T$ and since $\alpha (G)$ $\leq $ $\alpha (T)=|A|=|B|\leq \alpha (G)$,
we infer that $\alpha (G)=|A|=|B|$, i.e., $G$ has $A,B$ as stability systems
and therefore is $\alpha ^{+}$-stable, by the above Lemma \ref{lem3.4}. \rule%
{2mm}{2mm}\newline

The bipartite complement of the bipartite graph $G$ $=(A,B,E)$ is the graph
denoted by $\widetilde{G}=$ $(A,B,\widetilde{E})$, where $\widetilde{E}=$ $%
\{ab:a\in A,b\in B$ and $ab\notin E\}$.

\begin{lemma}
\label{lem3.7} If $G$ is a connected balanced bipartite graph, then the
following conditions are equivalent:

$\mathit{(i)}$ $\mu (G+e)=\mu (G)$, for any $e\in E(\overline{G})$;

$\mathit{(ii)}$ $\mu (G+e)=\mu (G)$, for any $e\in E(\widetilde{G})$.
\end{lemma}

\setlength {\parindent}{0.0cm}\textbf{Proof. }$\mathit{(i)}$ $\Rightarrow $\ 
$\mathit{(ii)}$ It is true, because $E(\widetilde{G})\subseteq E(\overline{G}%
)$.\setlength {\parindent}{3.45ex}

$\mathit{(ii)}$ $\Rightarrow $\ $\mathit{(i)}$ Since $G$ is balanced and $%
\mu (G+e)=\mu (G)$, for any $e\in E(\widetilde{G})$, we infer that $G$ has
no unmatched pair of vertices, and therefore $\mu (G+e)=\mu (G)$, for any $%
e\in E(\overline{G})$. \rule{2mm}{2mm}

\begin{theorem}
\label{th3.8} If $G$ is a connected bipartite graph, then the following
assertions are equivalent:

$\mathit{(i)}$ $G$ is $\alpha ^{+}$-stable;

$\mathit{(ii)}$ $G$ has a perfect matching;

$\mathit{(iii)}$ $G$ possesses two stability systems that partition its
vertex set;

$\mathit{(iv)}$ $\cap \{S:S$ is a stability system of $G\}=$ $\emptyset $;

$\mathit{(v)}$ $|A|=|B|$ and $\mu (G+e)=\mu (G)$, for any $e\in E(\overline{G%
})$;

$\mathit{(vi)}$ $|A|=|B|$ and $\mu (G+e)=\mu (G)$, for any $e\in E(%
\widetilde{G})$.
\end{theorem}

\setlength {\parindent}{0.0cm}\textbf{Proof. }$\mathit{(i)}$ $\Rightarrow $\ 
$\mathit{(ii)}$ If $G$ is $\alpha ^{+}$-stable, Proposition \ref{prop3.6}
ensures that $G$ contains an $\alpha ^{+}$-stable spanning tree $T$. By
Theorem \ref{th3.2}, $T$ has a perfect matching, which clearly is a perfect
one for $G$ itself, since it covers $V(T)=V(G)$, using edges contained only
in $E(T)\subseteq E(G)$.\setlength {\parindent}{3.45ex}

$\mathit{(ii)}$ $\Rightarrow $\ $\mathit{(iii)}$ Suppose $G$, with
bipartition $\{A,B\}$, has a perfect matching $M$. Clearly, $M$ uses edges
from $(A,B)$ and consequently, $|A|=|B|=\alpha (G)$. Hence, $A$ and $B$ are
two stability systems for $G$ that partition its vertex set.

$\mathit{(iii)}$ $\Rightarrow $\ $\mathit{(iv)}$ It is clear.

$\mathit{(iv)}$ $\Rightarrow $\ $\mathit{(i)}$ It is true according to
Theorem \ref{th3.1}.

$\mathit{(ii)}$ $\Rightarrow $\ $(v)$ If $G$ has a perfect matching, then $%
|A|$ $=$ $|B|$ and there is no unmatched pair of vertices in $G$. Therefore,
we have also $\mu (G+e)=\mu (G)$, for any $e\in E(\overline{G})$.

$(v)$ $\Rightarrow $\ $\mathit{(ii)}$ Suppose, on the contrary, that $G$ has
no perfect matching. Consequently, for an arbitrary maximum matching $M$ of $%
G$, there are at least two unmatched vertices, say $a\in A$ and $b\in B$.
Clearly, $e=ab\notin E(G)$, otherwise $M\cup \{e\}$ is a matching in $G$
with more edges than $M$. Therefore, $e\in E(\overline{G})$ and $M\cup \{e\}$
is a matching in $G+e$; hence it follows that $\mu (G+e)>\mu (G)$,
contradicting the assumption in $\mathit{(v)}$.

$\mathit{(v)}$ $\Leftrightarrow $ $\mathit{(vi)}$ It is true according to
the above Lemma \ref{lem3.7}. \rule{2mm}{2mm}\newline

An interesting consequence of Theorem \ref{th3.8} is the following result:

\begin{proposition}
\label{prop3.9}If $G$ is a connected bipartite graph of order at least two,
then $|\cap \{S:S$ is a stability system of $G\}|$ $\neq $ $1$.
\end{proposition}

\setlength {\parindent}{0.0cm}\textbf{Proof. }Suppose, on the contrary, that
there is some connected bipartite graph $G$ with at least two vertices, such
that $|\cap \{S:S$ is a stability system of $G\}|=1$. According to Theorems 
\ref{th3.1} and \ref{th3.8}, $G$ is $\alpha ^{+}$-stable and has $|\cap
\{S:S $ is a stability system of $G\}|=0$, in contradiction with our
assumption, and this completes the proof. \rule{2mm}{2mm}%
\setlength
{\parindent}{3.45ex}

\begin{corollary}
\label{cor3.10}If $G$ is connected, $|\cap \{S:S$ is a stability system of $%
G\}|=1$, and $\left| V\left( G\right) \right| \geq 2$, then $G$ is
non-bipartite.
\end{corollary}

Any $K_{n},n\geq 2$ is $\alpha ^{+}$-stable and has $|\cap \{S:S$ is a
stability system of $G\}|=0$. For $n\geq $ $3$, the graph $H$ having $%
V(H)=V(K_{n})\cup \{v\}$ and $E(H)=E(K_{n})\cup \{vw\}$, with $w\in V(K_{n})$%
, is also $\alpha ^{+}$-stable, but has $|\cap \{S:S$ is a stability system
of $G\}|=1$.

\begin{corollary}
\label{cor3.11}For any $n\geq $ $4$ there exist connected non-bipartite $%
\alpha ^{+}$-stable graphs $G_{1}$ and $G_{2}$ of order $n$, such that $%
|\cap \{S:S$ is a stability system of $G_{1}\}|=0$, and $|\cap \{S:S$ is a
stability system of $G_{2}\}|=$ $1$.
\end{corollary}

\section{Bistable bipartite graphs}

If $G=(A,B,E)$ has $A$ and $B$ as its exactly two stability systems, then $G$
is a bistable bipartite graph (see \cite{LM97}). As we shall see in the next
section, this special class of bipartite graphs plays a key role in
describing the structure of $\alpha ^{+}$-stable and $\alpha $-stable
bipartite graphs.

It is easy to see that a bistable bipartite graph is necessarily connected.
Clearly, no chordless path $P_{n},n\geq $ $3$, is bistable bipartite.
However, it is possible to ''join'' some kind of chordless paths to a
bistable bipartite graph in such a way that the result is also bistable
bipartite. More precisely, we have:

\begin{lemma}
\label{lem4.1}Let $H$ be a bistable bipartite graph and $G=H\cup P_{k},k\geq
2$, be bipartite. Then $G$ is bistable bipartite if and only if $k$ is an
even number and at least the pendant vertices of $P_{k}$ are joined to some
vertices of $H$.
\end{lemma}

\setlength {\parindent}{0.0cm}\textbf{Proof. }Let $H=(A,B,E),P_{k}=(X,Y,U)$
with $x_{1},x_{k}$ as the pendant vertices, $x_{1}\in X$, and $G=(A\cup
X,B\cup Y,W)$.\setlength {\parindent}{3.45ex}

''\textit{if}'' Since both $G$ and $H$ are bistable bipartite, they have an
even number of vertices, and therefore $k$ must be even. If one of $%
(x_{1},B),(x_{k},A)$, say the first, is empty, then the set $B\cup Y\cup
\{x_{1}\}-\{x_{2}\}$ is a third stability system of $G$, a contradiction,
because $G$ is bistable bipartite. Therefore, we may conclude that $k$ is
even and at least the pendant vertices of $P_{k}$ are joined to some
vertices of $H$.

''\textit{only if}'' Since $G$ is connected, we have $G$ $=(A\cup X,B\cup
Y,W)$ and $\{A\cup X,B\cup Y\}$ is its unique bipartition (such that any of
its edges joins a vertex from $A\cup X$ to a vertex in $B\cup Y$). Assume
that $G$ has a stability system $S$, and $S\neq $ $A\cup X,S\neq B\cup Y$.
Let us denote $S_{A}=S\cap A,S_{B}=S\cap B,S_{X}=S\cap X,S_{Y}=S\cap Y$.

\textit{Case 1}. Both $S_{A}$, $S_{B}$ are non-empty. Then $%
|S_{A}|+|S_{B}|<|A|,|S_{X}|+|S_{Y}|\leq |X|$, and hence the contradiction: 
\[
\alpha (G)=|S|=|S_{A}|+|S_{B}|+|S_{X}|+|S_{Y}|<|A|+|X|=|A\cup X|\leq \alpha
(G). 
\]

\textit{Case 2}. $S_{A}$ $=A,S_{B}=\emptyset ,S_{X}\neq \emptyset ,S_{Y}\neq
\emptyset $. Then $x_{k}\notin S$, and there exists $i\in \{1,...,k-1\}$
such that $x_{i},x_{i+1}\notin S$. Henceforth, we get that $%
|S_{X}|+|S_{Y}|<|X|$, and again the same contradiction: 
\[
\alpha (G)=|S|=|S_{A}|+|S_{X}|+|S_{Y}|<|A|+|X|=|A\cup X|\leq \alpha (G). 
\]
Analogously for the case $S_{A}=\emptyset ,S_{B}=B,S_{X}\neq \emptyset
,S_{Y}\neq \emptyset $.

Consequently, we may have only $S=A\cup X$ or $S=B\cup Y$, and clearly it
follows that $G$ is bistable bipartite. \rule{2mm}{2mm}\newline

A \textit{vertex cover} of a graph $G$ is a set $W\subseteq V(G)$ such that $%
W$ contains at least an endpoint of every edge of $G$.

\begin{lemma}
\label{lem4.2}A set $S$ of vertices is stable in $G$ if and only if $V(G)-S$
is a vertex cover of $G$, and the cardinality of a minimum vertex cover
equals $|V(G)|-\alpha (G)$.
\end{lemma}

\setlength {\parindent}{0.0cm}\textbf{Proof. }By definition, $S$ is stable
if and only if no edge of $G$ joins two of its vertices, and this implies
that any edge of $G$ has at least one of its endpoints in $V(G)-S$, i.e., $%
V(G)-S$ is a vertex cover of $G$. Conversely, if $V(G)-S$ is a vertex cover
in $G$, then no edge of $G$ has both endpoints in $S$, that is $S$ is stable
in $G$. Consequently, $S$ is a stability system of $G$ if and only if $%
V(G)-S $ is a minimum vertex cover of $G$, and clearly the cardinality of a
minimum vertex cover of $G$ equals $|V(G)|-\alpha (G)$. \rule{2mm}{2mm}%
\setlength
{\parindent}{3.45ex}\newline

A bipartite graph $H=(A,B,E)$ is said to be \textit{cover-irreducible} if it
is balanced and $A,B$ are its only minimum vertex covers, (see \cite{DM67}).
In \cite{LP77} a graph $G$ is defined as \textit{elementary} if the union of
all its perfect matchings forms a connected subgraph of $G$. The next
theorem extends the characterization of elementary bipartite graphs, due to
Hetyei (see \cite{H64}), Lovasz and Plummer (see \cite{LP77}).

\begin{theorem}
\label{th4.3}If $G=(A,B,E)$ is a bipartite graph with at least $4$ vertices,
then the following assertions are equivalent:

$\mathit{(i)}$ $G$ is cover-irreducible;

$\mathit{(ii)}$ $G$ is bistable;

$\mathit{(iii)}$ for any proper subset $X$ of $A$ or of $B$, $|N(x)|>|X|$
holds;

$\mathit{(iv)}$ $G$ is balanced and for any proper subset $X$ of $A$, $%
|N(x)|>|X|$ holds;

$\mathit{(v)}$ $G-a-b$ is $\alpha ^{+}$-stable, for any $a\in A$ and $b\in B$%
;

$\mathit{(vi)}$ for any $a\in A$ and $b\in B,G-a-b$ has a perfect matching;

$\mathit{(vii)}$ $G$ is connected and any of its edges lies in some perfect
matching of $G$;

$\mathit{(viii)}$ $G$ is elementary;

$\mathit{(ix)}$ $G$ can be written in the form $G=G_{0}\cup H_{1}\cup
H_{2}\cup ...\cup H_{k}$, where $G_{0}$ consists of two vertices and an edge
joining them and $H_{i}$ is an even path which joins two points of $%
G_{0}\cup H_{1}\cup ...\cup H_{i-1}$ in different color classes and has no
other point in common with $G_{0}\cup H_{1}\cup ...\cup H_{i-1}$.
\end{theorem}

\setlength {\parindent}{0.0cm}\textbf{Proof. }$\mathit{(i)}$ $%
\Leftrightarrow $ $\mathit{(ii)}$ By definition, $G$ is bistable if and only
if it has $A$ and $B$ as its only two stability systems, and according to
Lemma \ref{lem4.2}, if and only if $V(G)-A=B$ and $V(G)-B=A$ are its only
two minimum vertex covers, i.e., $G$ is cover-irreducible.%
\setlength
{\parindent}{3.45ex}

$\mathit{(ii)}$ $\Rightarrow $\ $\mathit{(iii)}$ If $A$ and $B$ are
stability systems in $G$, then clearly $\alpha (G)=|A|=|B|$. Suppose there
is some proper subset $X$ of $A$ such that $|N(x)|\leq |X|$. Consequently, $%
(X,B-N(x))=\emptyset $, and therefore $S=X\cup (B-N(x))$ is a stable set in $%
G$, which has $|S|=|X|+|B-N(x)|\geq |X|+|A-X|=\alpha (G)$. Hence, $S$ is a
third stability system of $G$, since it meets both $A$ and $B$, a
contradiction. Analogously, if $X\subset B$.

$\mathit{(iii)}$ $\Rightarrow $\ $\mathit{(iv)}$ Suppose that $|A|<|B|$.
Then, for any $b\in B,|B-\{b\}|\geq |A|$, and according to $\mathit{(iii)}$,
we must have $|A|\leq |B-\{b\}|<|N(B-\{b\})|\leq |A|$, a contradiction.
Analogously, the inverse inequality is also impossible, and hence, $|A|=|B|$
holds. Thus, $\mathit{(iii)}$ is true.

$\mathit{(iv)}$ $\Rightarrow $\ $\mathit{(v)}$ Let $a\in A,b\in B$ be and $%
H=G-a-b$. By Theorem \ref{th3.8}, it suffices to show that $\alpha
(H)=|A-\{a\}|=|B-\{b\}|$. Suppose, on the contrary, that $\alpha (H)=\alpha
(G)$; then there is a stable set $S$ in $H$, such that $\alpha
(H)=|S|>|A-\{a\}|$. Then both $S_{A}=$ $S\cap (A-\{a\})$ and $S_{B}=$ $S\cap
(B-\{b\})$ are non-empty. Since $N(S_{A})\subseteq B-\{b\}-S_{B}$, we obtain
the following contradiction: $|N(S_{A})|\leq
|B-\{b\}-S_{B}|<|S_{A}|<|N(S_{A})|$.

Therefore, $\alpha (H)=|A-\{a\}|=|B-\{b\}|$ must hold.

$\mathit{(v)}$ $\Rightarrow $\ $\mathit{(vi)}$ It is true, according to
Theorem \ref{th3.8}.

$\mathit{(vi)}$ $\Rightarrow $\ $\mathit{(vii)}$ Clearly, $G$ is connected,
since otherwise for $a,b$ in different color classes and different connected
components $G-a-b$ has no perfect matching, contradicting the assumption on $%
G-a-b$. Let $ab$ be an arbitrary edge of $G$ and $M$ be a perfect matching
in $G-a-b$, which exists according to hypothesis. Hence, $M\cup \{ab\}$ is a
perfect matching in $G$ containing $ab$.

$\mathit{(vii)}$ $\Rightarrow $\ $\mathit{(viii)}$ It is clear, since the
subgraph of $G$ formed by $\cup \{M:M$ is a perfect matching of $G\}$ is $G$
itself.

$\mathit{(viii)}$ $\Rightarrow $\ $\mathit{(ix)}$ It has been proved in \cite
{LP77}.

$\mathit{(ix)}$ $\Rightarrow $\ $\mathit{(ii)}$ We use induction on the
number $r$ of the even paths $H_{i}$. Clearly, $G_{0}\cup H_{1}$ is
bistable, since it is an even cycle. According to Lemma \ref{lem4.1}, if $%
G_{0}\cup H_{1}\cup ...\cup H_{r-1}$ is bistable and $H_{r}$ is an even
path, then by its construction, $G_{0}\cup H_{1}\cup ...\cup H_{r-1}\cup
H_{r}$ is also bistable. \rule{2mm}{2mm}

\begin{corollary}
\label{cor4.4}Any bistable bipartite graph having at least $4$ vertices is $%
\alpha $-stable.
\end{corollary}

\setlength {\parindent}{0.0cm}\textbf{Proof. }Let $G$ $=(A,B,E)$ be a
bistable bipartite graph. According to Theorem 4.3$\mathit{(v)}$, $G$ is $%
\alpha ^{+}$-stable. By Theorem \ref{th4.3}$\mathit{(iii)}$, $%
|N(v)|>|\{v\}|=1$ holds for any vertex $v$ of $G$, and this implies that
both $A$ and $B$ are $2$-dominating. Consequently, $G$ is also $\alpha ^{-}$%
--stable, by Theorem \ref{th2.2}. \rule{2mm}{2mm}%
\setlength
{\parindent}{3.45ex}\newline

Combining Corollary \ref{cor4.4} and Theorems \ref{th3.8} and \ref{th2.14},
we get the following result:

\begin{corollary}
\label{cor4.5}Any bistable bipartite graph $G$ with at least $4$ vertices
has perfect matchings and $\cap \{M:M$ is a perfect matching of $%
G\}=\emptyset $.
\end{corollary}

In the next proposition we define some special kind of ''bipartite''
substitution, which also preserves the ''bistable bipartite'' property.

\begin{proposition}
\label{prop4.6}Let $G_{i}=(A_{i},B_{i},E_{i}),i=1,...,p$, be disjoint
bistable bipartite graphs, $H=(X,Y,U)$ be a bipartite graph having $%
X=\{a_{1},...,a_{p}\},Y=\{b_{1},...,b_{p}\}$ and $G=(A,B,E)$ be defined as
follows: $A=A_{1}\cup ...\cup A_{p},B=B_{1}\cup ...\cup B_{p}$ and $%
E=E_{1}\cup ...\cup E_{p}\cup \{a_{i}b_{j}:$ for some $a_{i}\in
A_{i},b_{j}\in B_{j}$ if $x_{i}y_{j}\in U,1\leq i,j\leq p,i\neq j\}$. Then $H
$ is bistable if and only if $G$ has the same property.
\end{proposition}

\setlength {\parindent}{0.0cm}\textbf{Proof. }Let $H$ be with $X,Y$ as its
only two stability systems and suppose that $G$ has a third stability system 
$S$, different from $A$ and $B$. Then, there is some $k\in \{1,...,p\}$,
such that both $S\cap A_{k}\neq \emptyset $ and $S\cap B_{k}\neq \emptyset $%
. Since for any $i\in \{1,...,p\}$, both $S\cap A_{i}$ and $S\cap B_{i}$ are
stable in $G_{i}$, we obtain: $|S\cap A_{i}|+|S\cap B_{i}|\leq
|A_{i}|=|B_{i}|$, if $i\neq k$, and $|S\cap A_{k}|+|S\cap
B_{k}|<|A_{k}|=|B_{k}|$, and hence the contradiction: 
\[
\alpha (G)=|S|=|S\cap A_{1}|+|S\cap B_{1}|+...+|S\cap A_{p}|+|S\cap
B_{p}|<|A_{i}|+...+|A_{p}|=\alpha (G). 
\]

\setlength {\parindent}{3.45ex}Conversely, assume that $G$ is bistable, but $%
H$ has a third stability system $S$. Therefore, both $S_{X}=S\cap X$ and $%
S_{Y}=S\cap Y$ are non-empty. Since $(S_{X},S_{Y})=\emptyset $, it follows
that $(\cup \{$ $A_{i}:a_{i}\in S_{X}\},\cup \{B_{i}:b_{i}\in S_{Y}\})=$ $%
\emptyset $, and consequently, the set $(\cup \{A_{i}:a_{i}\in S_{X}\})\cup
(\cup \{B_{i}:b_{i}\in S_{Y}\})$ is a third stability system of $G$, which
contradicts our assumption on $G$. \rule{2mm}{2mm}

\begin{corollary}
\label{cor4.7}Let $G=(A,B,E)$ be a connected bipartite graph which possesses
a partition of its vertex set consisting of two simple cycles $C_{1},C_{2}$
and let denote $A_{i}=C_{i}\cap A$ and $B_{i}=C_{i}\cap B$, for $i=1,2$.
Then $G$ is bistable bipartite if and only if both $(A_{1},B_{2})\neq
\emptyset $ and $(A_{2},B_{1})\neq \emptyset $.
\end{corollary}

\setlength {\parindent}{0.0cm}\textbf{Proof. }Now, the graph $H$ from the
above Proposition \ref{prop4.6} must be isomorphic to a chordless cycle on $%
4 $ vertices, i.e., both $(A_{1},B_{2})\neq \emptyset $ and $%
(A_{2},B_{1})\neq \emptyset $ hold in $G$. \rule{2mm}{2mm}%
\setlength
{\parindent}{3.45ex}

\section{$\alpha $-Stable bipartite graphs}

A graph $G=(V,E)$ is called $\alpha $\textit{-stable} if it is both $\alpha
^{-}$-stable and $\alpha ^{+}$-stable, (see \cite{LM97}). In this section we
determine the structure of $\alpha $-stable and of $\alpha ^{+}$-stable
bipartite graphs in terms of bistable bipartite graphs. We also provide
constructions for larger $\alpha $-stable and $\alpha ^{+}$-stable bipartite
graphs from smaller ones.

A disconnected graph $G$, with components $H_{1},...,H_{p}$, is $\alpha $%
-stable if and only if the next assertions are valid:

$\mathit{(i)}$ each $H_{i},(i=1,...,p)$, is $\alpha $-stable;

$\mathit{(ii)}$ at most one of $H_{i},i=1,...,p$, has $|\cap \{S_{i}:S_{i}$
is a stability system in $H_{i}\}|=1$.

Combining Theorems \ref{th2.5} and \ref{th3.1} we obtain the following
result:

\begin{proposition}
\label{prop5.1}\cite{LM97} There is no connected chordal graph with at least
two vertices which is $\alpha $-stable.
\end{proposition}

In particular, no tree of order at least two is $\alpha $-stable. However
there exist connected $\alpha $-stable bipartite graphs (e.g., any even
chordless cycle). Moreover, Theorems \ref{th2.14} and \ref{th3.8} yield the
following characterization for the $\alpha $-stable bipartite graphs:

\begin{theorem}
\label{th5.2}If $G$ is a connected bipartite graph, then the following
assertions are equivalent:

$\mathit{(i)}$ $G$ is $\alpha $-stable;

$\mathit{(ii)}$ $G$ has perfect matchings and $\cap \{M:M$ is a perfect
matching of $G\}=\emptyset $;

$\mathit{(iii)}$ $|A|=|B|$ and $\mu (G-e)=\mu (G)=\mu (G+u)$, for any $e\in E
$ and $u\in E(\overline{G})$;

$\mathit{(iv)}$ $|A|=|B|$ and $\mu (G-e)=\mu (G)=\mu (G+u)$, for any $e\in E$
and $u\in E(\widetilde{G})$;

$\mathit{(v)}$ $\cap \{S:S$ is stable, $\left| S\right| =\alpha \left(
G\right) \}=\cap \{M:M$ is a matching, $\left| M\right| =\mu \left( G\right)
\}=$ $\emptyset $.
\end{theorem}

\begin{lemma}
\label{lem5.3}If graph $G=(A,B,E)$ is $\alpha $-stable and $S$ is a
stability system of $G$ meeting both $A$ and $B$, then the subgraph $%
H=G[(S\cap A)\cup (B-S)]$ is $\alpha $-stable.
\end{lemma}

\setlength {\parindent}{0.0cm}\textbf{Proof. }Since $S_{A}$ $=$ $S\cap A$
and $B-S_{B}$, (for $S_{B}=S\cap B$), are matched in any perfect matching of 
$G$, $H$ is $\alpha ^{+}$-stable, by Theorem \ref{th3.8}. We show that $H$
is also $\alpha ^{-}$--stable. Both $A$ and $B$ are $2$-dominating in $G$,
and therefore $H$ has at least $4$ vertices. $S_{A}$ is $2$-dominating,
because for any $b\in B-S_{B}$, we have $|N(b)\cap S_{A}|=|N(b)\cap S|\geq 2$%
. $S_{B}$ is also $2$-dominating, since for any $a\in B-S_{A}$, we have $%
|N(a)\cap S_{B}|=|N(a)\cap S|\geq 2$. Let $X$ be a stability system of $H$,
with $X_{A}=X\cap A=X\cap S_{A}\neq \emptyset $ and $X_{B}=X\cap B=X\cap
(B-S_{B})\neq $ $\emptyset $ . Then $S^{\prime }=X\cup S_{B}$ is clearly a
stability system of $G$, and hence, we have: $|N(a)\cap X|=|N(a)\cap
X_{B}|=|N(a)\cap S^{\prime }|\geq 2$, for any $a\in S_{A}-X_{A}$, and $%
|N(b)\cap X|=|N(b)\cap X_{A}|=|N(a)\cap S^{\prime }|\geq 2$, for any $b\in
B-S_{B}-X_{B}$, i.e., $X$ is $2$-dominating in $H$. Consequently, $H$ is $%
\alpha ^{-}$--stable, by Theorem \ref{th2.2}. \rule{2mm}{2mm}%
\setlength
{\parindent}{3.45ex}

\begin{lemma}
\label{lem5.4}If a bipartite graph has two perfect matchings $M_{1},M_{2}$,
then any of its vertices, from which are issuing edges contained in $%
M_{1},M_{2}$, respectively, belongs to some cycle that is alternating with
respect to at least one of $M_{1},M_{2}$.
\end{lemma}

\setlength {\parindent}{0.0cm}\textbf{Proof. }Let $G=(A,B,E)$ be a bipartite
graph having $M_{1}$, $M_{2}$ as perfect matchings, and let $a\in A$ be some
vertex from which are issuing two edges, say $ab_{1}$ and $ab_{2}$,
contained in $M_{1},M_{2}$, respectively. Henceforth, there exist $%
a_{2}b_{2}\in M_{1}$ and $a_{1}b_{1}\in M_{2}$. If $a_{1}=a_{2}$ or $%
a_{1}b_{2}\in E(G)$, we are done. Otherwise, we have $a_{1}b_{3}\in M_{1}$
and $a_{2}b_{4}\in M_{2}$. If $b_{3}=b_{4}$ or $a_{2}b_{3}\in E(G)$ or $%
a_{1}b_{4}\in E(G)$, we are done. If this does not happen, we continue with $%
a_{3}b_{3}\in M_{2}$ and $a_{4}b_{4}\in M_{1}$. If $a_{3}=a_{4}$ or $%
a_{4}b_{1}\in E(G)$ or $a_{4}b_{3}\in E(G)$ or $a_{3}b_{2}\in E(G)$ or $%
a_{3}b_{4}\in E(G)$, we are done. Otherwise, we continue in the same way.
Since $G$ is finite, the process must end, and a cycle, alternating with
respect to $M_{1}$ or $M_{2}$, is revealed. \rule{2mm}{2mm}%
\setlength
{\parindent}{3.45ex}

\begin{corollary}
\label{cor5.5}A bipartite graph $G$\ has two disjoint perfect matchings if
and only if it has a partition of its vertex set comprising of a family of
simple cycles.
\end{corollary}

\setlength {\parindent}{0.0cm}\textbf{Proof. }If $x$ is a vertex in $G$,
then we can build a cycle $C_{x}=(A_{x},B_{x},E_{x})$ using only edges
belonging alternately to the two perfect matchings $M_{1},M_{2}$, as is
shown in the proof of the above lemma. If $y$ is a vertex not contained in $%
C_{x}$, we use the same procedure to create a new cycle $C_{y}$, alternating
with respect to the both perfect matchings. This cycle does not use vertices
of the former cycle, since $A_{x}$ and $B_{x}$ are already matched by the
two perfect matchings and $C_{y}$ uses only edges of $M_{1}$ and $M_{2}$. In
this way, we get a partition of $V(G)$ consisting of vertex sets of pairwise
disjoint cycles.\setlength {\parindent}{3.45ex}

Conversely, let $\{C_{i}:1\leq i\leq k\}$ be a family of cycles of $G$, such
that their vertex sets $\{V(C_{i}):1\leq i\leq k\}$ form a partition for $%
V(G)$. Since $G$ is bipartite, all $C_{i}$ are even and each one has two
disjoint perfect matchings, say $M_{i,1}$ and $M_{i,2}$. Hence, $\cup
\{M_{i,1}:1\leq i\leq k\}$ and $\cup \{M_{i,2}:1\leq i\leq k\}$ are two
disjoint perfect matchings of $G$ itself. \rule{2mm}{2mm}\newline

The main result is as follows:

\begin{theorem}
\label{th5.6}If $G$ is a connected bipartite graph, then the following
assertions are equivalent:

$\mathit{(i)}$ $G$ is $\alpha $-stable;

$\mathit{(ii)}$ $G$ can be decomposed as $G$ $=G_{1}\cup ...\cup G_{k},k\geq
1$, such that each $G_{i},1\leq i\leq k$, is bistable bipartite and has at
least $4$ vertices;

$\mathit{(iii)}$ $G$ has perfect matchings and $\cap \{M:M$ is a perfect
matching of $G\}=\emptyset $;

$\mathit{(iv)}$ from each vertex of $G$ are issuing at least two edges,
which are contained in some perfect matchings of $G$;

$\mathit{(v)}$ any of its vertices belongs to some alternating cycle of $G$.
\end{theorem}

\setlength {\parindent}{0.0cm}\textbf{Proof. }$\mathit{(i)}$ $\Rightarrow $\ 
$\mathit{(ii)}$ If $G=(A,B,E)$ has $A$ and $B$ as its only two stability
systems, we are done. Otherwise, let $S$ be a stability system of $G$, such
that both $S_{A}=S\cap A$ and $S_{B}=S\cap B$ are non-empty. By Lemma \ref
{lem5.3}, the subgraphs $H_{1}$ $=$ $G[(S\cap A)\cup (B-S)]$ and $%
H_{2}=G[(A-S)\cup (S\cap B)]$ are $\alpha $-stable. If they both are
bistable, we are done. Otherwise, we continue with this decomposition
procedure, until all the subgraphs we obtain are bistable. After a finite
number of subpartitions, we obtain a decomposition of $G$ as $G=G_{1}\cup
G_{2}\cup ...\cup G_{k},k\geq 1$, such that each $G_{i},1\leq i\leq k$, is
bistable and $\alpha $-stable.\setlength
{\parindent}{3.45ex}

$\mathit{(ii)}$ $\Rightarrow $\ $\mathit{(iii)}$ By assertion $\mathit{(vi)}$
of Theorem \ref{th4.3}, each $G_{i},1\leq i\leq k$, has a perfect matching.
Hence $G$ owns a perfect matching, too. According to Corollary \ref{cor4.5},
it follows that $\cap \{M_{i}:M_{i}$ is a perfect matching of $%
G_{i}\}=\emptyset $ for every $1\leq i\leq k$. Therefore, $\cap \{M:M$ is a
perfect matching of $G\}\subseteq \cap \{\bigcup\limits_{i=1}^{k}$ $%
M_{i}:M_{i}$ is a perfect matching of $G_{i}\}=$ $\emptyset $, i.e., $\cap
\{M:M$ is a perfect matching of $G\}=\emptyset $ must hold.

$\mathit{(iii)}$ $\Rightarrow $\ $\mathit{(iv)}$ Suppose the conclusion is
not true. Then, there is a vertex $v$ in $G$, so that only one edge, say $vw$%
, is contained in some perfect matching of $G$; such an edge must exist,
because $G$ has perfect matchings. Moreover, since $v$ is matched with a
vertex by each perfect matching, we infer that $vw$ is contained in all
perfect matchings of $G$, in contradiction with the hypothesis on $G$.
Therefore, $\mathit{(iv)}$ is valid.

$\mathit{(iv)}$ $\Rightarrow $\ $\mathit{(v)}$ It is true, according Lemma 
\ref{lem5.4}.

$\mathit{(v)}$ $\Rightarrow $\ $\mathit{(i)}$ Since any vertex of $G=(A,B,E)$
is contained in some alternating cycle $C_{x}$ of $G$, we infer, according
to Theorem \ref{th3.8}, that $G$ is $\alpha ^{+}$-stable. It is easy to see
that $A$ and $B$ are $2$-dominating in $G$. Let $S$ be a stability system of 
$G$, with both $S_{A}=S\cap A$ and $S_{B}=S\cap B$ non-empty. Suppose, on
the contrary, that $S$ is not $2$-dominating in $G$, i.e., there is some
vertex $a\in A$, so that $|N(a)\cap S|=|\{b\}|=1$. Let $C_{a}$ be an
alternating cycle containing $a$. Since $A-S_{A}$ is matched with $S_{B}$ by
any perfect matching of $G$, the edge $ab$ belongs to $C_{a}$ and the second
neighbor of $a$ on $C_{a}$, say $c$, must be in $B-S_{B}$, otherwise we
obtain that $2=|\{b,c\}|\leq |N(a)\cap S_{B}|=|N(a)\cap S|$, in
contradiction with the assumption on $a$. Hence, because $%
(S_{A},S_{B})=\emptyset $, the cycle $C_{a}$ must use at least two edges
from $(A-S_{A},B-S_{B})$, which are contained in no perfect matching of $G$.
Therefore, $C_{a}$ can not be alternating in $G$, contradicting the choice
of $C_{a}$. Consequently, each stability system of $G$ is $2$-dominating,
and by Theorem \ref{th2.2}, $G$ is also $\alpha ^{-}$--stable. \rule%
{2mm}{2mm}\newline

As a consequence, we get the following result:

\begin{proposition}
\label{prop5.7}A connected bipartite graph $G$ is $\alpha ^{+}$-stable if
and only if it admits a decomposition as $G=G_{1}\cup ...\cup G_{k}$, where
all $G_{i}$ are bistable bipartite.
\end{proposition}

\setlength {\parindent}{0.0cm}\textbf{Proof. }Let $M_{0}=\cap \{M:M$ is a
perfect matching in $G\}$, $H_{0}$ be the subgraph of $G$ spanned by the
vertices matched by the edges contained in $M_{0}$, and $H_{1}=G-H_{0}$.%
\setlength {\parindent}{3.45ex}

Clearly, $H_{1}$ has $\cap \{M:M$ is a perfect matching in $H_{1}\}$ $=$ $%
\emptyset $, while $H_{0}$ is either empty or a disjoint union of $K_{2}$.
According to Theorem \ref{th5.6}, any connected component of $H_{1}$ has a
decomposition in bistable bipartite graphs. Therefore, $G$ admits a
decomposition as $G=G_{1}\cup ...\cup G_{k}$, all $G_{i}$ being bistable.

Conversely, if $G=G_{1}\cup ...\cup G_{k}$, and all $G_{i}$ are bistable,
then each $G_{i}$ has at least a perfect matching $M_{i}$, and $\cup
\{M_{i}:M_{i}$ is a perfect matching in $G_{i}\}$ is a perfect matching in $%
G $. Consequently, by Theorem \ref{th3.8}, $G$ is $\alpha ^{+}$-stable. \rule%
{2mm}{2mm}

\begin{proposition}
\label{prop5.8}Let $H=(A,B,E)$ be a connected, bipartite and $\alpha $%
-stable graph and $G=H\cup P_{k}$ , for $k\geq $ $3$, be connected. Then $G$
is $\alpha $-stable if and only if $k$ is even and for any stability system $%
S$ of $H$, at least one of the pendant vertices of $P_{k}$ is adjacent to
some vertex of $S$.
\end{proposition}

\setlength {\parindent}{0.0cm}\textbf{Proof. }Let $P_{k}=(X,Y,U)$ and $%
x_{1},x_{k}$ be its pendant vertices.\setlength {\parindent}{3.45ex}

Suppose $G$ is $\alpha $-stable. Hence, $|V(G)|$ is even, and consequently, $%
k=2p$. Assume $x_{1}\in X$ and $x_{k}\in Y$. By Theorem \ref{th3.8}, we have 
$\alpha (G)=|A\cup X|=|B\cup Y|=\alpha (H)+\alpha (P_{k})$. Hence, $A\cup X$
and $B\cup Y$ are $2$-dominating stability systems in $G$, and therefore $%
(x_{1},B)$ and $(x_{k},A)$ are both non-empty. Let $S$ be another stability
system of $G$. Since clearly $\alpha (G)=|S|=|S\cap V(H)|+|S\cap
V(P_{k})|\leq \alpha (H)+\alpha (P_{k})=\alpha (G)$, we get that $S\cap V(H)$
is a stability system in $H$ and $S\cap V(P_{k})$ is a stability system in $%
P_{k}$. If suppose $x_{1},x_{k}\in S$, there is $i\in \{1,...,k-3\}$ such
that $x_{i}\in S\cap X$ and $x_{i+3}\in S\cap Y$, and consequently we get $%
|N(x_{i+2})\cap S|=|\{x_{i+3}\}|=$ $1$, contradicting the fact that $S$ is $%
2 $-dominating in $G$. Hence, either $(x_{1},S)\neq \emptyset $ or $%
(x_{k},S)\neq $ $\emptyset $.

Conversely, by Theorem \ref{th3.8}, $G$ is $\alpha ^{+}$-stable, since it
has at least a perfect matching. In addition, we have that $\alpha
(G)=|A\cup X|=|B\cup Y|=\alpha (H)+\alpha (P_{k})$, and because $(x_{1},B)$
and $(x_{k},A)$ are both non-empty, we infer that $A\cup X$ and $B\cup Y$
are $2$-dominating stability systems in $G$. If $S$ is a another stability
system of $G$, as above, we get that $S\cap V(H)$ and $S\cap V(P_{k})$ are
stability systems in $H$ and $P_{k}$, respectively. According to the
hypothesis, one of $x_{1}$ and $x_{k}$, say $x_{k}$, is joined to some
vertex $t$ of $S\cap V(H)$. Consequently, we get that $S\cap V(P_{k})=X$.
Now, since $S\cap V(H)$ is $2$-dominating in $H$, and any $z\in Y-\{x_{k}\}$
has two neighbors in $X\subset S$, and $N(x_{k})\cap S=\{t,x_{k}-1\}$, we
may assert that $S$ is $2$-dominating in $G$. Consequently, $G$ is also $%
\alpha ^{-}$--stable, and finally $G$ is $\alpha $-stable. \rule{2mm}{2mm}

\begin{proposition}
\label{prop5.9}If $G_{i},1\leq i\leq k,k\geq $ $1$ are bipartite and $%
G=G_{1}\cup ...\cup G_{k}$ is connected, then:

$\mathit{(i)}$ $G$ is $\alpha ^{+}$-stable whenever all $G_{i},1\leq i\leq k$%
, are $\alpha ^{+}$-stable;

$\mathit{(ii)}$ $G$ is $\alpha $-stable whenever all $G_{i},1\leq i\leq k$,
are $\alpha $-stable.
\end{proposition}

\setlength {\parindent}{0.0cm}\textbf{Proof. }$\mathit{(i)}$ Since all $%
G_{i},1\leq i\leq k$, by Theorem \ref{th3.8}, have a perfect matching, $G$
itself has a perfect matching, and therefore it is $\alpha ^{+}$-stable.%
\setlength {\parindent}{3.45ex}

$\mathit{(ii)}$ If all $G_{i}=(A_{i},B_{i},E_{i}),1\leq i\leq k$, are $%
\alpha $-stable, then according to the part $\mathit{(i)}$, $G$ $=(A,B,E)$
is $\alpha ^{+}$-stable and 
\[
\alpha (G)=|A|=|\cup \{A_{i}:1\leq i\leq k\}|=|B|=|\{B_{i}:1\leq i\leq
k\}|=\alpha (G_{1})+...+\alpha (G_{k}). 
\]

Clearly, $A$ and $B$ are $2$-dominating in $G$, and if $S$ is another
stability system of $G$, we have $\alpha (G)$ $=$ $|S|$ $=$ $|S\cap
V(G_{1})|+...+|S\cap V(G_{k})|\leq \alpha (G_{1})+...+\alpha (G_{k})=\alpha
(G)$, and this implies that $|S\cap V(G_{i})|=\alpha (G_{i})$ for any $i\in
\{1,...,k\}$, i.e., $S\cap V(G_{i})$ is a stability system in $G_{i},i\in
\{1,...,k\}$. Hence $S\cap V(G_{i})$ is $2$-dominating in $G_{i}$, for any $%
i\in \{1,...,k\}$, and therefore $S=\cup \{S\cap V(G_{i}):1\leq i\leq k\}$
is $2$-dominating in $G$. Consequently, $G$ is also $\alpha ^{-}$--stable,
by Theorem \ref{th2.2}. So, we may conclude that $G$ is $\alpha $-stable. 
\rule{2mm}{2mm}

\section{Conclusions}

In this paper we investigated the $\alpha ^{-}$--stability and $\alpha ^{+}$%
-stability of connected chordal graphs and of connected bipartite graphs,
characterizing the $\alpha ^{-}$--stable chordal graphs, $\alpha ^{-}$%
--stable bipartite, $\alpha ^{+}$-stable bipartite and $\alpha $-stable
bipartite graphs. These findings generalize some previously known results
for trees. We present new facts on strong unique independence graphs, and
determine the structure of $\alpha ^{+}$-stable bipartite graphs and of the $%
\alpha $-stable bipartite graphs. Several operations preserving these
structures are also considered.

\section{Acknowledgments}

We are indebted to Professor M.C. Golumbic for carefully reading and
commenting a part of this paper. We also gratefully thank Professor J. Topp
for pointing out several inaccuracies in the first version of the article,
and for his remark that Proposition 2.1 is true not only for bipartite
graphs, but for general graphs as well.


\begin{thebibliography}{99}
\bibitem{BHNS83}  D. Bauer, F. Harary, J. Nieminen and C. L. Suffel, \emph{%
Domination alteration sets in graphs}, Discrete Mathematics \textbf{47}
(1983) 153-161.

\bibitem{B73}  C. Berge, \emph{Graphs and Hypergraphs}, North Holland,
Amsterdam (1973).

\bibitem{BCD88}  R. C. Brigham, P. Z. Chinn and R. D. Dutton, \emph{Vertex
domination-critical graphs}, Networks \textbf{18} (1988) 173-179.

\bibitem{DM67}  A.L. Dulmage and N.S. Mendelsohn, \emph{Graphs and Matrices}%
, in Graph Theory and Theoretical Physics, ed. by F. Harary, Academic Press,
1967.

\bibitem{DB88}  R. D. Dutton and R. C. Brigham, \emph{An extremal problem
for edge domination insensitive graphs}, Discrete Applied Mathematics 
\textbf{20} (1988) 113-125.

\bibitem{G80}  M.C. Golumbic, \emph{Algorithmic Graph Theory and Perfect
Graphs}, Academic Press, New York, 1980

\bibitem{F93}  O. Favaron, \emph{A note on the irredundance number after
vertex deletion}, Discrete Mathematics \textbf{121} (1993) 51-54.

\bibitem{FJ84}  F. Fink and M. S. Jacobson, $n$\emph{-Domination in graphs},
Graph Theory with Applications to Algorithms and Computer Science, Wiley
(1985) 283-300.

\bibitem{GHR93}  G. Gunther, B. Hartnell, and D. F. Rall, \emph{Graphs whose
vertex independence number is unaffected by single edge addition or deletion}%
, Discrete Applied Mathematics \textbf{46} (1993) 167-172.

\bibitem{HLBD90}  T. W. Haynes, L. M. Lawson, R. C. Brigham and R. D.
Dutton, \emph{Changing and unchanging of the graphical invariants: minimum
and maximum degree, maximum clique size, node independence number and edge
independence number}, Congressus Numerantium \textbf{72} (1990) 239-252.

\bibitem{HOS96}  M. A. Henning, O. R. Oellermann and H.C. Swart, \emph{Local
edge domination critical graphs}, Discrete Mathematics \textbf{161} (1996)
175-184.

\bibitem{H64}  G. Hetyei, \emph{2 x 1-es teglalapokkal lefedheto idomokrol},
Pecsi Tanarkepzo Foisk. Tud. Kozl. (1964) 351-368.

\bibitem{HS85}  G. Hopkins and W. Staton, \emph{Graphs with unique maximum
independent sets}, Discrete Mathematics \textbf{57} (1985) 245-251.

\bibitem{K31}  D. K\"{o}nig, \emph{Graphen und Matrizen}, Mat. Fiz. Lapok 
\textbf{38} (1931) 116-119.

\bibitem{LM97}  V. E. Levit and E. Mandrescu, \emph{On }$\alpha $\emph{%
-stable graphs}, Congressus Numerantium \textbf{124} (1997) 33-46.

\bibitem{LP77}  L. Lovasz and M. D. Plummer, \emph{On minimal elementary
bipartite graphs}, Journal of Combinatorial Theory, Series A\textbf{\ 23}
(1977) 127-138.

\bibitem{M96}  S. Monson, \emph{The effects of vertex deletion and edge
deletion on clique partition number}, Ars Combinatoria \textbf{42} (1996)
89-96.

\bibitem{STV94}  W. Siemes, J. Topp and L. Volkmann, \emph{On unique
independent sets in graphs}, Discrete Mathematics \textbf{131} (1994)
279-285.

\bibitem{SB83}  D.P.Sumner and P.Blitch, \emph{Domination critical graphs},
Journal of Combinatorial Theory, Series B \textbf{34} (1983) 65-76.

\bibitem{S90}  D. P. Sumner, \emph{Critical concepts in domination},
Discrete Mathematics \textbf{86} (1990) 33-46.

\bibitem{TV98}  J. Topp and P. D. Vestergaard, $\alpha _{k}$\emph{- and }$%
\beta _{k}$\emph{-stable graphs}, manuscript.
\end{thebibliography}
\end{document}